\documentclass[11pt,reqno]{amsart}
\usepackage{enumitem}
\usepackage{xcolor}
\newlength{\basicwidth}
\newlength{\shortbasicwidth}
\newlength{\basicheight}
\numberwithin{equation}{section}
\newtheorem{REM}{Remark}

\begin{document}

\newcommand{\Cr}{\color{red}}
\newcommand{\Cb}{\color{blue}}
\newcommand{\Label}{\label}

\begin{center}
{{\Large{Sturm Theorem and a Refinement of Vietoris' Inequality\\[1.4ex]
for Cosine Polynomials}}}
\end{center}

\vspace{0.6cm}
\begin{center}
HORST ALZER$^a$  \quad{and}
\quad{MAN KAM KWONG$^b$}
\footnote{The research of this author is supported by the Hong Kong Government GRF Grant PolyU 5003/12P and the Hong Kong Polytechnic University Grants G-UC22 and G-UA10}
\end{center}

\vspace{0.6cm}
\begin{center}
$^a$ Morsbacher Str. 10, 51545 Waldbr\"ol, Germany\\
\emph{e-mail:} {\tt{h.alzer@gmx.de}}
\end{center}

\vspace{0.5cm}
\begin{center}
$^b$ Department of Mathematics, The Hong Kong Polytechnic University,\\
Hunghom, Hong Kong\\
\emph{e-mail:} {\tt{mankwong@poly.edu.hk}}
\end{center}

\vspace{0.8cm}
\begin{quote}
{\small 
{\bf{Abstract.}} 
In a recent work \cite{ak}, the authors established the following refinement of the
well-known 1958 result of Vietoris on nonnegative cosine polynomials.

\par\vspace*{6mm}\par

{\Cb
\noindent
{\bf Proposition. }
\em
Let

\par\vspace*{-11mm}\par
$$
T_n(x)=\sum_{k=0}^n b_k\cos(kx)
$$
with 
$$  \qquad  \qquad  b_{2k}=b_{2k+1}=\frac{1}{4^k}{\binom {2k}{k} }  \qquad (k\geq 0).   \eqno (*)  $$
The inequalities
$$
T_n(x)\geq c_0 + c_1 x + c_2 x^2 >0
\qquad{(c_k\in\mathbf{R}, k=0,1,2)}
$$
hold for all  $n\geq 1$ and  $x\in (0,\pi)$ if and only if
$$
c_0=\pi^2 c_2, 
\quad{c_1=-2\pi c_2,}
\quad{0<c_2\leq \alpha},
$$
where
$$
\alpha=\min_{0\leq t<\pi} 
\frac{T_6(t)}{(t-\pi)^2}=0.12290... .
$$
}

\par\vspace*{4mm}\par

\noindent
In four places of the proof, use was made of the classical Sturm Theorem
on determining the number of real roots of an algebraic polynomial
in a given interval.  Although absolutely rigorous, the
Sturm procedure involves lengthy technical computations carried
out with the help of the software MAPLE 13. This article
supplements \cite{ak} by providing such details which
were omitted in the latter.
}
\end{quote}

\vspace{0.8cm}

\vspace{0.1cm}

\newpage

\section{Introduction}

Please refer to the Abstract for the motivation of writing this
article.

Let 
\par\vspace*{-13mm}\par
\begin{equation}  X_0(y)=\sum_{k=0}^{n} c_ky^k   \end{equation}
be an AP (algebraic polynomial) with real coefficients $ c_k $,
$k=1,2,\cdots,n$, and $ (\alpha ,\beta )\subset\mathbb{R} $ be a given subinterval
of the real line. The celebrated Sturm Theorem furnishes a rigorous procedure
to determine the number of real roots (multiple roots are counted
once) of $ X_0(y) $ in $ (\alpha ,\beta ) $. Details of the
Sturm theorem are given in
van der Waerden \cite[p. 248 ]{W}; see also Kwong \cite{Kw}.

The steps of the SP (Sturm procedure) are summarized as follows. 
Suppose that $ \alpha  $ and $ \beta  $ are not roots of $ X(y) $, which is true in most 
applications.
\begin{enumerate}
\item[1.] Compute $ X_1(y)=X_0'(y) $, the derivative of $ X_0(y) $.
\item[2.] Compute the Sturm sequence of polynomials $ X_2(y),\,X_3(y),\,\cdots $
  using the Euclidean algorithm:
  Each $ X_i(y) $ is the negative of
  the remainder when $ X_{i-2} $ is divided by $ X_{i-1} $. In the algorithm
  used by MAPLE, each $ X_i $ is normalized (by dividing by a positive constant)
  so that the leading coefficient is $ \pm 1 $.
\item[3.] Count the number of sign changes in
  each of the two sequences $ \left\{ X_i(\alpha )\right\}  $ and $ \left\{ X_i(\beta )\right\}  $.
\item[4.] The difference of these two numbers is the number of real roots of $ X(y) $ in 
  $ (\alpha ,\beta ) $.
\end{enumerate}

The procedure can be easily modified if $ \alpha  $ and{/}or $ \beta  $ are roots of
$ X(y) $.

MAPLE provides a simple command to automate the procedure.
It displays only the desirable output from step 4, and hides all the 
intermediate, less relevant data from steps 1 to 3.

Every CP (cosine polynomial) $ \sum b_k\,\cos(kx) $
can be rewritten as an AP of the variable $ y=\cos(x) $ while
every SP (sine polynomial) $ \sum a_k\,\sin(kx) $
is the product of $ \sin(x) $ and an AP of $ y $.
This simple observation has enabled the authors to successfully exploit
the SP in helping (often only a relatively small portion of the entire
proof is based on the SP) to obtain new results on NN (nonnegative) 
TP (trigonometric polynomials), some of which
have been described in \cite{Kw}.

Intensive computations are involved in the practical
implementation of the SP. In our study, we have used
the software MAPLE 13 to carry out these computations.
As we have adamantly pointed out in \cite{Kw}, proofs 
based on the SP are absolutely theoretically rigorous.
However, some researchers consider the details of such computations to be 
too technical and uninteresting to be included in the presentation 
of the results in a theoretical research
article. With this in mind when we wrote \cite{ak}, we have omitted
all such details, and referred interested readers to this article.

The authors do not intend to publish this article in a regular journal.
It will be archived in arXiv and permanently available in the internet.

\par\vspace*{-2mm}\par
\section{Use of SP in \cite{ak}}

The following applications of the SP appear in the same order as in \cite{ak}. 

\begin{enumerate}[leftmargin=0.3in,itemsep=10pt]
\setlength{\parskip}{1.5ex}
\setlength{\parindent}{0pt}

\item {\Cb {\bf{Lemma 1.}}
\par\vspace*{-10mm}\par
\begin{equation}
\min_{0\leq t<\pi}\frac{T_6(x)}{(x-\pi)^2}
=0.1229... 
\end{equation}
\em attained at $ x=0.726656896\cdots $.
}

\bigskip
The proof makes use of the following functions
\begin{equation}  \eta(x)=10x^6+6x^5-12x^4-\frac{11}{2}x^3+\frac{29}{8}x^2+\frac{11}{8} x+\frac{9}{16}  \end{equation}
\begin{equation}  \quad{\mu(x)=\eta(x)\eta''(x)-\frac{1}{2}\eta'(x)^2}.  \end{equation}
and
\par\vspace*{-10mm}\par
\begin{equation}  \nu(x) = (1-x^2)^{3}\mu(x)^2 -4\,(0.1229)\, x^2 \eta(x)^3  \end{equation}
and states the claim that $\mu$ and
$\nu$ have no zeros in $[0.65,0.95]$.

Using MAPLE, we obtain the explicit expressions
{\footnotesize 
\begin{eqnarray*}
 \mu(x) &=& 1200\,{x}^{10}+1200\,{x}^{9}-1890\,{x}^{8}-1854\,{x}^{7}+938\,{x}^{6}+
987\,{x}^{5}+  \\[1.2ex]
&& \qquad  {\frac {615}{8}}\,{x}^{4}-{\frac {835}{8}}\,{x}^{3}-{
\frac {1659}{16}}\,{x}^{2}-{\frac {297}{16}}\,x+{\frac {401}{128}}
\end{eqnarray*}
}
and

\par\vspace*{-7mm}\par

{\footnotesize 
\begin{eqnarray*}
\nu (x) &=& -1440000\,{x}^{26}-2880000\,{x}^{25}+7416000\,{x}^{24}+17625600\,{x}^{
23} -14981700\,{x}^{22} - \\[1.5ex]
  && \qquad  47224920\,{x}^{21}+ {\frac {62018162}{5}}\,{x}^{
20}+{\frac {1799880978}{25}}\,{x}^{19}+{\frac {922239083}{250}}\,{x}^{
18} -\\[1.2ex]
  && \qquad  {\frac {84199650417}{1250}}\,{x}^{17}- 
  {\frac {92661304227}{5000}}
\,{x}^{16}+{\frac {190171138621}{5000}}\,{x}^{15}+ \\[1.2ex]
  && \qquad  {\frac {762930005877
}{40000}}\,{x}^{14}-{\frac {107540100801}{10000}}\,{x}^{13}-{\frac {
810190932293}{80000}}\,{x}^{12} -\\[1.2ex]
  && \qquad  {\frac {33113449049}{80000}}\,{x}^{11}
+{\frac {880096473123}{320000}}\,{x}^{10}+{\frac {103309397713}{80000}
}\,{x}^{9} -\\[1.2ex]
  && \qquad  {\frac {236121121411}{1280000}}\,{x}^{8}-{\frac {
425398321107}{1280000}}\,{x}^{7}-{\frac {840336154901}{10240000}}\,{x}
^{6}+ \\[1.2ex]
  && \qquad  {\frac {19241962691}{1280000}}\,{x}^{5}+{\frac {164343549777}{
10240000}}\,{x}^{4}+{\frac {18143739883}{5120000}}\,{x}^{3}- \\[1.2ex]
  && \qquad  {\frac {
428328477}{1280000}}\,{x}^{2}-{\frac {119097}{1024}}\,x+{\frac {160801
}{16384}} \, .
\end{eqnarray*}
}

Suppose that in a MAPLE session, the variables ``{\tt mu}'' and ``{\tt nu}''
have been assigned the polynomials $ \mu (x) $ and $ \nu (x) $, respectively. 
Applying the SP to the two polynomials is accomplished simply by the
commands

\begin{verbatim}
     sturm(mu, x, 65/100, 95/100); sturm(nu, x, 65/100, 95/100);
\end{verbatim}

The output are both 0, meaning that $ \mu (x) $ and $ \nu (x) $ have
no zeros in $ [0.65,0.95] $.
That is all we are interested to know.

Note that in the above commands, we have to use $ 65/100 $ and $ 95/100 $
in place of the decimal forms $ 0.65 $ and $ 0.95 $, respectively.
If the decimal forms are used, MAPLE will perform the SP
using floating point arithmetic instead of exact arithmetic. In the
floating point mode, rounding errors can and very often
do lead to erroneous results.

MAPLE actually makes use of the Sturm sequence in its internal computation, 
without
displaying them explicitly. If one is curious, one can see them using
the command

\begin{verbatim}
     sturmseq(mu, x);
\end{verbatim}

Just for the sake of illustration, 
we list below the first few members of the Sturm sequence 
of $ \mu (x) $, starting with $ \mu _1(x)=\mu '(x) $.
{\footnotesize 
\par\vspace*{-4mm}\par
$$  \mu _1 = {x}^{9}+{\frac {9}{10}}\,{x}^{8}-{\frac {63}{50}}\,{x}^{7}-{\frac { 2163}{2000}}\,{x}^{6}+{\frac {469}{1000}}\,{x}^{5}+{\frac {329}{800}} \,{x}^{4}+{\frac {41}{1600}}\,{x}^{3}-{\frac {167}{6400}}\,{x}^{2}-{ \frac {553}{32000}}\,x-{\frac {99}{64000}}  $$
$$  \mu _2 = {x}^{8}+5/6\,{x}^{7}-{\frac {25249}{24300}}\,{x}^{6}-{\frac {2429}{ 2700}}\,{x}^{5}+{\frac {43}{6480}}\,{x}^{4}+{\frac {6091}{38880}}\,{x} ^{3}+{\frac {473}{2880}}\,{x}^{2}+{\frac {1951}{64800}}\,x-{\frac { 10619}{1555200}} \hspace*{4mm}  $$
$$  \mu _3 = {x}^{7}+{\frac {164171}{403140}}\,{x}^{6}-{\frac {253857}{134380}}\,{x }^{5}-{\frac {24703}{26876}}\,{x}^{4}+{\frac {86929}{161256}}\,{x}^{3} +{\frac {104433}{430016}}\,{x}^{2}+{\frac {24223}{537520}}\,x+{\frac { 50933}{12900480}} \hspace*{2mm}  $$
$$  \hspace*{1mm}  \mu _4 = -{x}^{6}-{\frac {3617836265259}{2968641792581}}\,{x}^{5}+{\frac { 617789163720}{2968641792581}}\,{x}^{4}+{\frac {5544754493585}{ 11874567170324}}\,{x}^{3}-{\frac {1101526123005}{47498268681296}}\,{x} ^{2} -  $$
$$  {\frac {244252761843}{23749134340648}}\,x+{\frac {1195017616507}{ 94996537362592}} \hspace*{50mm}  $$
}

The two sequences of numbers mentioned in step 3 of the SP are
\par\vspace*{-5mm}\par
{\footnotesize 
$$  \left\{ \mu _i(0.65)\right\}  = \left\{ 0.006,\,\, 0.01,\,\, 0.01,\,\, -0.01,\,\, -0.05,\,\, 0.01,\,\, 0.01,\,\, -0.4,\,\, -0.6,\,\, -0.9,\,\, -1 \right\}   $$
$$  \left\{ \mu _i(0.95)\right\}  = \left\{ 0.2,\,\, 0.2,\,\, 0.09,\,\, -0.4,\,\, -1.1,\,\, 0.9,\,\, -0.1,\,\, -1.3,\,\, -1.3,\,\, -1.2,\,\, -1 \right\}   $$
}

\par\vspace*{-6mm}\par
There are 3 changes of sign in the first sequence, after the third, fifth, and
seventh terms, respectively.
There are also 3 changes of sign in the second sequence, after the third, 
fifth, and sixth terms, respectively. Hence, the outcome of step 4 is 0.

We are lucky that MAPLE hides all such gory details from us.

\item {\Cb {\bf{Lemma 7.}} \emph{Let
$$
C_n(x)=\sum_{k=0}^n (-1)^k \, b_k \cos(kx).
$$
{\rm(}where $ b_k $ are defined by {\rm(*)} in the Abstract.{\rm)}
If $2\leq n\leq 21$ $(n\neq 6)$ and $x\in (5\pi/8,\pi)$, then}
\begin{equation}  \frac{820}{33}\Bigl(1-\cos\frac{x}{10}\Bigr) \leq C_n(x).  \Label{l7}  \end{equation}
}

\par\vspace*{-4mm}\par
\begin{proof}
 We
set $y=x/10$ and
$$
P_n(y)=C_n(10y)-\frac{820}{33}(1-\cos(y)).
$$
Letting $Y=\cos(y)$. Then $P_n(y)$ is
 an algebraic polynomial in $Y$.
We denote this polynomial by $P_n^*(Y)$, where
 $Y\in [\cos(\pi/10), \cos(\pi/16)]
=[0.951...,0.980...]$. Applying Sturm's
theorem gives that $P_n^*$ has
no zero on $[0.951, 0.981]$ and satisfies $P_n^*(0.97)>0$. It follows that $P_n$
is positive on $[\pi/16, \pi/10]$. This implies that (\ref{l7}) holds. 
\end{proof}

Let us use the case $ n=2 $ as an illustration. We have
\begin{equation}  P_2(y) = 1 - \cos(10y) + \frac{1}{2} \,\cos(20y) - \frac{820}{33} \,(1-\cos(y)) 1  \end{equation}
Assume that the MAPLE variable ``{\tt P}'' has already been assigned this cosine polynomial.
The command

\par\vspace*{-2mm}\par

\begin{verbatim}
     X := subs(cos(y)=Y,expand(P));
\end{verbatim}
\par\vspace*{-3mm}\par
produces the AP $ P_3^*(Y) $ mentioned in the proof
and assigns it to the variable ``{\tt X}''.

\par\vspace*{-5mm}\par
{\footnotesize 
$$  X=262144\,{Y}^{20}-1310720\,{Y}^{18}+2785280\,{Y}^{16}-3276800\,{Y}^{14} +2329600\,{Y}^{12}-1025536\,{Y}^{10}+  $$
$$  \hspace*{10mm} 275840\,{Y}^{8}-43360\,{Y}^{6}+ 3700\,{Y}^{4}-150\,{Y}^{2}+{\frac {820}{33}}\,Y-{\frac {1475}{66}}  $$
}

\par\vspace*{-6mm}\par
Our initial goal is to affirm that $ X $ is NN in the interval 
$Y\in [\cos(\pi/10), \cos(\pi/16)]$. However, the MAPLE command ``{\tt sturm}''
can only handle intervals where the endpoints $ \alpha  $ and $ \beta  $ are rational
numbers. To overcome that obstacle, we show that $ X $ is NN in
the slightly larger interval $[0.951, 0.981]$ using the command
\begin{verbatim}
     sturm(X, Y, 951/1000, 981/1000);
\end{verbatim}

The same process has to be repeated for the other cases $ 2<n\leq 21 $
($n\neq 6$). The following MAPLE snippet can be used for that purpose.

\begin{verbatim}
     b := proc (n) local n1; 
        if n = 1 or n = 0 then 1 
        else n1 := floor((1/2)*n)+1; 
           factorial(2*n1-3)/(2^(2*n1-3)*factorial(n1-1)*factorial(n1-2)) 
        end if 
     end proc;

     C := n -> 1 + add((-1)^k*b(k)*cos(k*x), k = 1 .. n);

     for n from 2 to 21 do
        P := C(n) - 820/33 * (1-cos((1/10)*x));
        X := subs(cos(y) = Y, expand(subs(x = 10*y, expand(P))));
        print((n,sturm(X, Y, 951/1000, 981/1000)));
     end do:
\end{verbatim}

The first six lines define a procedure to generate the coefficients of the
Vietoris polynomial: the command ``{\tt b(n)}'' produces 
$ b_n $ as defined in ($\ast$). 
The next line defines a procedure to generate $ C_n(x) $ defined in 
Lemma 7. The remaining lines constitute a do loop to apply the SP to
verify (\ref{l7}), for $ n $ from 2 to 21.

The output consists of pairs of numbers, 
``{\tt n}'' and the number of roots of $ X $ in the interval in question. The latter is 0
except for $ n=6 $. That is why this case is excluded in Lemma 7.

\begin{REM} \em
The statement (fourth line from the bottom)
used to generate $ X $ from $ P $ is different from that given 
earlier. This is because the MAPLE command ``{\tt expand}'' only knows
how to expand $ \cos(kx) $ when $ k<100 $. The earlier statement
will fail if $ n\geq 10 $. The modified statement expands
the expression before making the substitution $ y=x/10 $. Then after
the substitution another expansion is performed to obtain $ X $.
\end{REM}

\begin{REM} \em
As $ n $ increases, the degrees of the corresponding $ P $ and $ X $ increase
and the time needed by the SP increases significantly. The SP is 
therefore not recommended for dealing with TP of high degrees. 

In fact, the SP is seldom adequate when establishing general results. At best,
it plays an assistant role used to take care of a few exceptional
cases involving lower degree polynomials.
\end{REM}

\begin{REM} \em
The symbol ``{\tt :}'' is used to end the do loop in the
last line to suppress printing of the intermediate computation results, namely,
the various $ P $ and $ X $ in the iterations. If one is curious to see
what these results are, one can use the regular statement terminator 
``{\tt;}'' instead. When $ n $ is large, the expression $ X $ is rather lengthy.
\end{REM}

\item {\Cb {\bf{Lemma 8.}} \emph{Let}
$$  \Delta(x)=\sum_{k=0}^{21} (-1)^k (b_k -b_{22}) \cos(kx)-\frac{820}{33}\Bigl(1-\cos\frac{x}{10}\Bigr).  $$

\par\vspace*{-8mm}\par
\begin{tabbing}\hspace*{4ex} \= \hspace*{5ex} \= \hspace*{20ex} \= \\
    \> If  \>  \hspace*{0.16mm} $5\pi/8\leq x\leq 2.68$,  \>  then \, $\Delta(x)> 0.29$; \\
    \> if  \>  \hspace*{1.1mm} $2.68\leq x\leq 2.83$,    \>  then \, $\Delta(x)> 0.46$; \\
    \> if  \>  \hspace*{1.1mm} $2.83\leq x\leq 2.908$,   \>  then \, $\Delta(x)> 0.64$; \\
    \> if  \>  $2.908\leq x\leq 2.970$,  \>  then \, $\Delta(x)> 0.90$; \\
    \> if  \>  $2.970\leq x\leq 3.021$,  \>  then \, $\Delta(x)> 1.32$; \\
    \> if  \>  $3.021\leq x\leq 3.051$,  \>  then \, $\Delta(x)> 1.78$. \\
\end{tabbing}
}

\par\vspace*{-4mm}\par
\begin{proof}
Let $5\pi/8\leq x\leq 2.68$. We have
$\cos(\pi/16)=0.980...$ and $\cos(0.268)=0.964...$. 
The function $\Delta -0.29$ is an AP
in $Y=\cos(x/10)$. An application of
Sturm's theorem shows that this function
is positive on $[0.964,0.981]$. This leads
to $\Delta(x)>0.29$ for $x\in [5\pi/8, 2.68]$. 
Using the same method of proof we obtain the other 
estimates for $\Delta(x)$.
\end{proof}

Essentially this Lemma reveals good lower bounds of $ \Delta (x) $ in progressively
smaller subintervals of $ [5\pi /8,\pi ] $.  Numerical optimization methods can
be used to estimate $ \min\left\{ \Delta (x)\right\}  $ in any given interval,
but that cannot be taken as a rigorous proof. 
Especially when the objective function is oscillatory in the interval,
as is the case of $ \Delta (x) $ on $ [5\pi /8,2.68] $,
a numerical algorithm may erroneously 
yield one of the larger local minimum instead of 
the global one. The SP provides a theoretically rigorous lower bound.

\item As usual, we adopt the notation
$$
(a)_0=1, \quad{(a)_n=\prod_{k=0}^{n-1}  (a+k)}=\frac{\Gamma(a+n)}{\Gamma(a)}
\quad{(n\geq 1)},
$$
and define
\par\vspace*{-9mm}\par
\begin{equation}
d_{2k} = d_{2k+1} = \frac{(69/100 )_k}{k! }
\qquad{(k=0,1,2,...)} .
\end{equation}

{\Cb{\bf{Lemma 11.}} \emph{Let}
\begin{equation}
I(x)= 
\sum_{k=0}^{21}\Bigl(b_k-\frac{b_{22}}{d_{22}} d_k\Bigr)\cos(kx).
\end{equation}
\emph{If $0<x\leq 0.1$, then} $I(x)>1.5$.
}

\begin{proof}
With $Y=\cos(x)\in[\cos(0.1),1]\subset[0.995,1]$, $(I(x)-1.5)$  
is an AP in $Y$. The SP reveals 
that it is positive in the interval.
It follows that
$I(x)>1.5$ for $x\in[0,0.1]$. 
\end{proof}

Here we are dealing with a single CP, $ (I(x)-1 5) $.
Again to avoid an irrational endpoint, we use the slightly larger
interval $ [995/1000,1] $ for the corresponding AP. Suppose that $ I(x) $
has been assigned to the variable ``{\tt Ix}'',
the MAPLE command to use is

\begin{verbatim}
     sturm(subs(cos(x)=Y, expand(Ix)), Y, 995/1000, 1);
\end{verbatim}

\end{enumerate}

\end{document}